\documentstyle[12pt]{amsart}
\begin{document}
\title{Biharmonic hypersurfaces in a Riemannian manifold with non-positive Ricci curvature}
\author{Nobumitsu Nakauchi and Hajime URAKAWA${}^{\ast}$ }
\title[Biharmonic hypersurfaces]
{Biharmonic hypersurfaces in a Riemannian manifold with non-positive Ricci curvature}

\thanks{$^{\ast}$\,\,Supported by the Grant-in-Aid for the Scientific Research, 
(C), No. 21540207, 
Japan Society for the Promotion of Science.}
\dedicatory{}
\maketitle
\begin{abstract} 
In this paper, we show that, for a biharmonic hypersurface $(M,g)$ 
of a Riemannian manifold $(N,h)$
of non-positive Ricci curvature, if 
$\int_M\vert H\vert^2 v_g<\infty$, 
where $H$ is the mean curvature of $(M,g)$ in $(N,h)$, 
then $(M,g)$ is minimal in $(N,h)$. 
Thus, for a counter example $(M,g)$ 
in the case of hypersurfaces to 
the generalized Chen's conjecture (cf. Sect.1), 
it holds that $\int_M\vert H\vert^2 v_g=\infty$. 
\end{abstract}
\numberwithin{equation}{section}
\theoremstyle{plain}
\newtheorem{df}{Definition}[section]
\newtheorem{th}[df]{Theorem}
\newtheorem{prop}[df]{Proposition}
\newtheorem{lem}[df]{Lemma}
\newtheorem{cor}[df]{Corollary}
\newtheorem{rem}[df]{Remark}
\section{Introduction and statement of results.} 
In this paper, we consider an isometric immersion 
$\varphi:\,(M,h)\rightarrow (N,h)$, 
of a Riemannian manifold $(M,g)$ of dimension $m$, 
into another Riemannian manifold 
$(N,h)$ of dimension $n=m+1$. We have
$$
\nabla^N_{\varphi_{\ast}X}\varphi_{\ast}Y=\varphi_{\ast}(\nabla_XY)+k(X,Y)\xi,
$$
for vector fields $X$ and $Y$ on $M$, 
where $\nabla$, $\nabla^N$  
are the Levi-Civita connections of 
$(M,g)$ and $(N,h)$, respectively, 
$\xi$ is the unit normal vector field along $\varphi$, 
and $k$ is the second fundamental form. 
Let 
$A: T_xM\rightarrow T_xM$ $(x\in M)$ be the shape operator defined by 
$g(AX, Y)=k(X,Y)$, $(X,Y\in T_xM)$, and $H$, the mean curvature 
defined by  $H:=\frac{1}{m}{\rm Tr}_g(A)$. 
Then, let us recall the following 
{\bf B.Y. Chen's conjectrure} (cf. \cite{C}, \cite{C2}): 
\vskip0.3cm\par
{\it Let $\varphi:\, (M,g)\rightarrow ({\mathbb R}^n, g_0)$ be an isometric immersion into the standard Euclidean space. If $\varphi$ is biharmonic 
(see Sect. 2), then, it is minimal. }
\vskip0.6cm\par
This conjecture is still open up to now, and let us recall also the following 
{\bf generalized B.Y. Chen's conjecture} (cf. \cite{C}, \cite{CMP}): 
\vskip0.3cm\par
{\it Let 
$\varphi:\,(M,g)\rightarrow (N,h)$ be an isometric immersion, and 
the sectional curvature of $(N,h)$ is non-positive. 
If $\varphi$ is biharmonic, then, it is minimal. }
\vskip0.6cm\par
Oniciuc (\cite{On}) and Ou (\cite{Ou2}) showed 
this is true if $H$ is constant. 
\vskip0.6cm\par
In this paper, we show 
\begin{th}
Assume that $(M,g)$ is complete 
and 
the Ricci tensor ${\rm Ric}^N$ of $(N,h)$ satisfies that 
\begin{equation}
{\rm Ric}^N(\xi,\xi)\leq \vert A\vert^2.
\end{equation}
If $\varphi:\,(M,g)\rightarrow (N,h)$ is biharmonic (cf. Sect. 2) and satisfies that 
\begin{equation}
\int_MH^2\,v_g<\infty,
\end{equation}
then, $\varphi$ has constant mean curvature,  i.e., $H$ is constant. 
\end{th} 
\vskip0.6cm\par
As a direct corollary, we have 
\begin{cor}
Assume that $(M,g)$ is a complete Riemannian manifold 
of dimension $m$ 
and $(N,h)$ is a Riemannian manifold of dimension 
$m+1$ whose Ricci curvature 
is non-positive. 
If
 an isometric immersion $\varphi:\,  (M,g)\rightarrow (N,h)$ is 
 biharmonic
and satisfies that $\int_MH^2\,v_g<\infty$, then, 
$\varphi$ is minimal. 
\end{cor}
\vskip0.6cm\par
By our Corollary 1.2, if there would exist a counter example (cf. \cite{Ou2}) 
in the case $\dim N=\dim M+1$, 
then 
it must hold that 
\begin{equation}
\int_MH^2\,v_g=\infty,
\end{equation}
which imposes the strong condition on the behaviour of the boundary of $M$ at infinity. 
Indeed, $(1.3)$ implies that either
$H$ is unbounded on $M$, 
or 
it holds that 
$H^2\geq C$ on an open subset $\Omega$ of $M$ with infinite volume, for some constant $C>0$. 
\vskip0.6cm\par
{\bf Acknowledgement.} 
\quad We express our thanks to the referee(s) who informed relevant references and gave useful comments to us.  
\vskip0.6cm\par
\section{Preliminaries.}
In this section, we prepare general 
materials about 
harmonic maps and biharmonic maps of a 
complete Riemannian manifold into another Riemannian manifold (cf. \cite{EL}). 
\par
Let $(M,g)$ be an $m$-dimensional complete Riemannian manifold, and 
the target space $(N,h)$ is 
an $n$-dimensional Riemannian manifold. 
For every $C^{\infty}$ map $\varphi$ of 
$M$ into $N$, and relatively compact domain $\Omega$ in $M$, 
the {\it energy functional} 
on the space $C^{\infty}(M,N)$ of all $C^{\infty}$ maps 
of $M$ into $N$ is defined by 
$$E_{\Omega}(\varphi)=\frac12\int_{\Omega}\vert d\varphi\vert^2\,v_g,$$
and 
for a $C^{\infty}$ one parameter deformation 
$\varphi_t\in C^{\infty}(M,N)$ $(-\epsilon<t<\epsilon)$ 
of  $\varphi$ with $\varphi_0=\varphi$, 
the variation vector field $V$ along $\varphi$ is defined by 
$V=\frac{d}{dt}\big\vert_{t=0} 
\varphi_t$.  
Let  $\Gamma_{\Omega}(\varphi^{-1}TN)$ be 
the space 
of $C^{\infty}$ sections of the induced bundle 
$\varphi^{-1}TN$ of the tangent bundle $TN$ by $\varphi$ 
whose supports are contained in $\Omega$. 
For $V\in \Gamma_{\Omega}(\varphi^{-1}TN)$ and 
its one-parameter deformation $\varphi_t$, 
the {\it first variation formula} is given by 
$$
\frac{d}{dt}
\bigg\vert _{t=0}E_{\Omega}(\varphi_t)
=-\int_{\Omega}
\langle 
\tau(\varphi),V \rangle \,v_g. 
$$ 
The {\it tension field} $\tau(\varphi)$ is 
defined globally on $M$ by 
\begin{equation}
\tau(\varphi)=\sum_{i=1}^mB(\varphi)(e_i,e_i),
\end{equation}
where
$$B(\varphi)(X,Y)=\nabla^N_{\varphi_{\ast}(X)}\varphi_{\ast}(Y)
-\varphi_{\ast}(\nabla_XY)$$
for $X,Y\in {\frak X}(M)$. 
Then, a $C^{\infty}$ map 
$\varphi:(M,g)\rightarrow (N,h)$ is {\it harmonic} 
if $\tau(\varphi)=0$. 
For a harmonic map $\varphi:\,(M,g)\rightarrow (N,h)$, 
the {\it second variation formula} of the energy functional 
$E_{\Omega}(\varphi)$ is 
$$
\frac{d^2}{dt^2}\bigg\vert_{t=0}
E_{\Omega}(\varphi_t)=\int_{\Omega}\langle J(V),V\rangle\,v_g
$$
where 
\begin{align*}
J(V)&:=\overline{\Delta}V-{\mathcal R}(V),\\
\overline{\Delta}V&
:={\overline{\nabla}}^{\ast}\,{\overline{\nabla}}V
=
-\sum_{i=1}^m
\{
{\overline{\nabla}}_{e_i}({\overline{\nabla}}_{e_i}V)
-{\overline{\nabla}}_{\nabla_{e_i}e_i}V
\},\\
{\mathcal R}(V)&:=\sum_{i=1}^m
R^N(V,\varphi_{\ast}(e_i))\varphi_{\ast}(e_i).
\end{align*}
Here, ${\overline{\nabla}}$ is the induced connection 
on the induced bundle $\varphi^{-1}TN$, and 
$R^N$ is the curvature tensor of $(N,h)$ 
given by 
$R^N(U,V)W=
[\nabla^N_{\,\,U},\nabla^N_{\,\,V}]W-\nabla^N_{\,\,[U,V]}W$ 
$(U,V,W\in {\frak X}(N)$).  
\par
The {\it bienergy functional} is 
defined by 
\begin{equation*}
E_{2,\Omega}(\varphi)=
\frac12\int_{\Omega}\vert\tau(\varphi)\vert^2\,v_g, 
\end{equation*}
and the {\it first variation formula} of the bienergy 
is given (cf. \cite{J}) by 
\begin{equation*}
\frac{d}{dt}\bigg\vert _{t=0}E_{2,\Omega}(\varphi_t)
=-\int_{\Omega}\langle \tau_2(\varphi),V\rangle\,v_g
\end{equation*}
where the {\it bitension field} $\tau_2(\varphi)$ is 
defined globally on $M$ by 
\begin{equation}
\tau_2(\varphi)=J(\tau(\varphi))
=\overline{\Delta}\tau(\varphi)-{\mathcal R}(\tau(\varphi)),
\end{equation}
and a $C^{\infty}$ map
$\varphi:(M,g)\rightarrow (N,h)$ is called to be 
{\it biharmonic} if 
\begin{equation}
\tau_2(\varphi)=0.
\end{equation}
\vskip0.6cm\par
\section{Some Lemma for the Schr\"odinger type equation}
In this section, we prepare some simple lemma 
of the Schr\"odinger type equation of 
the Laplacian 
$\Delta_g$ on an $m$-dimensional non-compact complete Riemannian manifold $(M,g)$ defined by 
\begin{equation}
\Delta_gf:=\sum_{i=1}^m
e_i(e_if)-\nabla_{e_i}e_if\qquad (f\in C^{\infty}(M)),
\end{equation}
where $\{e_i\}_{i=1}^m$ is a locally defined orthonormal frame field on $(M,g)$. 
\begin{lem}
Assume that $(M,g)$ is a complete non-compact 
Riemannian manifold, and $L$ is a non-negative 
smooth
function on $M$. 
Then, every smooth $L^2$ function $f$ on $M$ 
satisfying the Schr\"odinger type equation
\begin{equation}
\Delta_gf=L\,f\qquad (\text{on}\,\,M)
\end{equation} 
must be a constant. 
\end{lem}
\begin{pf} 
Take any point $x_0$ in $M$, and for every 
$r>0$, let us consider the following cut-off function 
$\eta$ on $M$: 
\begin{equation}
\left\{
\begin{aligned}
&0\leq\eta(x)\leq 1 \quad (x\in M),\\
&\eta(x)=1\quad (x\in B_r(x_0)),\\
&\eta(x)=0\quad (x\not\in B_{2r}(x_0)),\\
&\vert\nabla\eta\vert\leq\frac{2}{r}\quad (\text{on}\,\,M),
\end{aligned}
\right.
\end{equation}
where $B_r(x_0)=\{x\in M:\,d(x,x_0)<r\}$, and 
$d$ is the distance of $(M,g)$. 
Multiply $\eta^2\,f$ on $(3.2)$, and 
integrale it over $M$, we have 
\begin{equation}
\int_M(\eta^2\,f)\,\Delta_gf\,v_g=\int_ML\,\eta^2\,f^2\,v_g. 
\end{equation}
By the integration by part for the left hand side, 
we have 
\begin{equation}
\int_M(\eta^2\,f)\,\Delta_gf\,v_g
=-\int_Mg(\nabla(\eta^2\,f),\nabla f)\,v_g. 
\end{equation}
Here, we have 
\begin{align}
g(\nabla(\eta^2\,f),\nabla f)
&=2\eta\,f\,g(\nabla\eta,\nabla f)+\eta^2\,g(\nabla f,\nabla f)\nonumber\\
&=2\,\eta\,f\,\langle \nabla \eta,\nabla f\rangle
+\eta^2\,\vert \,\nabla f\,\vert^2,
\end{align}
where we use 
$\langle\,\cdot\,,\,\cdot\,\rangle$ and $\vert\,\cdot\,\vert$ 
instead of $g(\,\cdot\,,\,\cdot\,)$ and 
$g(u,u)=\vert u\vert^2$ $(u\in T_xM)$, for simplicity. 
Substitute $(3.6)$ into $(3.5)$, the right hand side of $(3.5)$ is equal to 
\begin{align}
RHS \,\text{of $(3.5)$}&=
-\int_M2\,\eta f\,\langle\nabla\eta,\nabla f\rangle\,v_g
-\int_M\eta^2\,\vert\nabla f\vert^2\,v_g\nonumber\\
&=-2\int_M\langle f\,\nabla\eta,\eta\,\nabla f\,\rangle\,v_g
-\int_M\eta^2\,\vert\nabla f\vert^2\,v_g.
\end{align}
\par
Here, applying Young's inequality:
for every $\epsilon>0$, and 
every vectors $X$ and $Y$ at each point of $M$,
\begin{equation}
\pm 2\,\langle X,Y\rangle
\leq\epsilon\vert X\vert^2+\frac{1}{\epsilon}\vert Y\vert^2, 
\end{equation}
 to the first term of $(3.7)$, we have 
\begin{align}
RHS\,\text{of $(3.7)$}
&\leq \epsilon\int_M\vert\eta\,\nabla f\vert^2\,v_g
+\frac{1}{\epsilon}\int_M\vert f\,\nabla \eta\vert^2\,v_g
-\int_M\eta^2\,\vert\nabla f\vert^2\,v_g\nonumber\\
&=-(1-\epsilon)\int_M\eta^2\,\vert\nabla f\vert^2\,v_g
+\frac{1}{\epsilon}\int_Mf^2\,\vert\nabla\eta\vert^2\,v_g.
\end{align}
Thus, by $(3.5)$ and $(3.9)$, we obtain
\begin{equation}
\int_ML\,\eta^2\,f^2\,v_g
+(1-\epsilon)\int_M\eta^2\,\vert\nabla f\vert^2\,v_g
\leq
\frac{1}{\epsilon}\int_Mf^2\,\vert\nabla\eta\vert^2\,v_g.
\end{equation}
\par
Now, puttig $\epsilon=\frac12$, $(3.10)$ implies that 
\begin{equation}
\int_ML\,\eta^2\,f^2\,v_g
+\frac12\int_M\eta^2\,\vert\nabla f\vert^2\,v_g
\leq
2\int_Mf^2\,\vert\nabla\eta\vert^2\,v_g.
\end{equation}
Since $\eta=1$ on $B_r(x_0)$ 
and $\vert\nabla\eta\vert\leq\frac{2}{r}$, and $L\geq 0$ on $M$, we have 
\begin{equation}
0\leq
\int_{B_r(x_0)}L\,f^2\,v_g
+\frac12\int_{B_r(x_0)}\vert\nabla f\vert^2\,v_g
\leq\frac{8}{r^2}\int_Mf^2\,v_g.
\end{equation} 
Since $(M,g)$ is non-compact and complete, $r$ can tend to infinity, and $B_r(x_0)$ goes to $M$. Then we have 
\begin{equation}
0\leq 
\int_ML\,f^2\,v_g+\frac12\int_M\vert\nabla f\vert^2\,v_g\leq 0
\end{equation}
since $\int_Mf^2\,v_g<\infty$. Thus, we have $L\,f^2=0$ and $\vert \nabla f\vert=0$ (on $M$) which implies that $f$ is a constant. 
\end{pf}
\vskip0.6cm\par
\section{Biharmonic isometric immersions.}
In this section, we consider 
a hypersurface $M$ of an $(m+1)$-dimensional Riemannian manifold $(N,h)$. Recently, Y-L. Ou showed (cf. 
\cite{Ou1}) 
\begin{th}
Let $\varphi:\,(M,g)\rightarrow (N,h)$ 
be an isometric immersion of an 
$m$-dimensional Riemannian manifold 
$(M,g)$ into another $(m+1)$-dimensional 
Riemannian manifold $(N,h)$ 
with the mean curvature vector field 
$\eta=H\,\xi$, 
where $\xi$ is the unit normal 
vector field along $\varphi$. Then, 
$\varphi$ is biharmonic if and only if the following 
equations hold:
\begin{equation}
\left\{
\begin{aligned}
&\Delta_gH-H\,\vert A\vert^2+H\,{\rm Ric}^N(\xi,\xi)=0,
\\
&2\,A\,(\nabla H)+\frac{m}{2}\,\nabla(H^2)-2\,H\,({\rm Ric}^N(\xi))^T=0,
\end{aligned}
\right.
\end{equation}
where 
${\rm Ric}^N:\,T_yN\rightarrow T_yN$ is the Ricci transform which is defined by 
$h({\rm Ric}^N(Z),W)={\rm Ric}^N(Z,W)$ 
$(Z,W\in T_yN)$, $(\cdot)^T$ is the tangential  component 
corresponding to the decomposition 
of $T_{\varphi(x)}N=\varphi_{\ast}(T_xM)\oplus 
{\mathbb R}\xi_x$ $(x\in M)$, 
and $\nabla f$ is the gradient vector field of 
$f\in C^{\infty}(M)$ on $(M,g)$, respectively.  
\end{th}
 \vskip0.6cm\par
 Due to Theorem 4.1 and Lemma 3.1, we can show immediately our Theorem 1.1. 
 \vskip0.3cm\par
({\it Proof of Theorem 1.1.})
\par \quad 
 Let us denote by $L:=\vert A\vert^2-{\rm Ric}^N(\xi,\xi)$ 
 which is a smooth non-negative function on $M$ due to our assumption. Then, the first equation is reduced to 
 the following Schr\"odinger type equation:
 \begin{equation}
 \Delta_gf=L\,f,
 \end{equation}
where $f:=H$ is a smooth $L^2$ function on $M$ by the assumption $(1.2)$. 
 \par
 Assume that $M$ is compact. 
 In this case,  by (4.2) and the integration by part, we have 
 \begin{equation}
 0\leq
 \int_ML\,f^2\,v_g=\int_Mf\,(\Delta_gf)\,v_g
 =-\int_Mg(\nabla f,\nabla f)\,v_g\leq 0,
 \end{equation}
 which implies that $\int_Mg(\nabla f,\nabla f)\,v_g=0$, 
 that is, $f$ is constant. 
 \par
 Assume that $M$ is non-compact. 
 In this case, we can apply Lemma 3.1 to (4.2). 
 Then, we have that $f=H$ is a constant.  
 \qed
 \vskip0.6cm\par
 ({\it Proof of Corollary 1.2.}) 
 \par
 Assume that ${\rm Ric}^N$ is non-positive. 
 Since $L=\vert A\vert^2-{\rm Ric}^N(\xi,\xi)$ is non-negative, $H$ is constant due to Theorem 1.1. 
 Then, due to $(4.1)$, 
 we have that 
 $H\,L=0$ and 
 $H\,({\rm Ric}^N(\xi))^T=0$. 
 If $H\not=0$, then $L=0$, i.e., 
 \begin{equation}
 {\rm Ric}^N(\xi,\xi)=\vert A\vert^2. 
 \end{equation}
By our assumption, ${\rm Ric}^N(\xi,\xi)\leq 0$, and 
the right hand side of $(4.4)$ is non-negative, so we have 
$\vert A\vert^2=0$, i.e., $A\equiv0$. 
This contradicts $H\not=0$. We have $H=0$. 
\qed 
\vskip2cm\par       

\vskip0.6cm\par
Graduate School of Science and Engineering, 
\par Yamaguchi University,
Yamaguchi, 753-8512, Japan.
\par
{\it E-mail address}: nakauchi@@yamaguchi-u.ac.jp
\vskip0.8cm\par
Division of Mathematics, Graduate School of Information Sciences, 
\par Tohoku University, 
Aoba 6-3-09, Sendai, 980-8579, Japan.
\vskip0.1cm\par
{\it Current Address}: 
Institute for International Education, 
\par Tohoku University,
Kawauchi 41, Sendai, 980-8576, Japan.
\par
{\it E-mail address}: urakawa@@math.is.tohoku.ac.jp
\end{document}